# Numerical Solution of Nonlinear Wave-Like Equations by Reduced Differential Transform Method

Murat Gubes[1+], Yildiray Keskin[2] and Galip Oturanc[2]

[1] Department of Mathematics, Kamil Ozdag Science Faculty of Karamanoglu Mehmetbey University, 70100, Campus, Karaman, Turkey

[2] Department of Mathematics, Science Faculty of Selcuk University, 42100, Campus-Konya, Turkey

## Abstract

This paper is presented to give numerical solutions of some cases of nonlinear wave-like equations with variable coefficients by using Reduced Differential Transform Method (RDTM). RDTM can be applied most of the physical, engineering, biological and etc. models as an alternative to obtain reliable and fastest converge, efficient approximations. Hence, our obtained results showed that RDTM is a very simple method and has a quite accuracy.

## 1. Introduction

Many physical problems can be described by mathematical models that involve partial differential equations. Large varieties of physical, chemical and biological phenomena are governed by the partial differential equations. A mathematical model is a simplified description of physical reality expressed in mathematical terms. Additionally, nonlinear partial differential equations are central to research in many fields such as Hydrodynamics, Engineering, Quantum field theory, Optics, Plasma physics etc. They mostly do not have exact solutions and therefore they are approximated using numerical schemes.

By applying the Adomian Decomposition Method, M. Ghoreishi solved some types of nonlinear wave-like equation [4], V.G. Gupta and S. Gupta worked out by using Homotopy Perturbation Transform Method these types of equation tool [5], furthermore, A. Aslanov [13], F. Yin and et al [12] and A. Atangana and et al [11] researched for solving nonlinear heat and wave-like equation by using Homotopy Perturbation, Variational Iteration and Homotopy Decomposition Methods respectively. Moreover, various techniques, such as homotopy analysis, perturbations, decompositions, iterations, differential and

---

[+] Corresponding author. Tel.: +90-338-226 2151 (3796); fax: +90-338-226 2150.
  *E-mail address*: mgubes@kmu.edu.tr and mrtgbs@gmail.com

laplace transformation techniques have been used to handle similar types of these wave-like and also heat-like problems numerically and analytically as in references [8-10], [13-17].

The fundamental motivation of the present study is the extension of a recently developed technique which is called Reduced Differential Transform Method (RDTM) to tackle some of nonlinear wave-like equations as the following form

$$u_{tt} = \sum_{i,j=1}^{n} F_{1ij}(X,t,u) \frac{\partial^{k+m}}{\partial x_i^k \partial x_j^m} F_{2ij}(u_{x_i}, u_{x_j})$$
$$+ \sum_{i=1}^{n} G_{1i}(X,t,u) \frac{\partial^p}{\partial x_i^p} G_{2i}(u_{x_i}) + H(X,t,u) + S(X,t) \tag{1}$$

with initial conditions

$$u(X,0) = a_0(X), \qquad u_t(X,0) = a_1(X). \tag{2}$$

Here, $X = (x_1, x_2, ..., x_n)$ and $F_{1ij}, G_{1i}$ are nonlinear functions of $X, t, u$. $F_{2ij}, G_{2i}$ are nonlinear functions of derivatives of $x_i$ and $x_j$ respectively. Also $H, S$ are nonlinear functions and $k, m, p$ are integers. These kind of equations describe the evolution of stochastic systems for example, erratic motions of small particles that are immersed in fluids, fluctuations of the intensity of laser light, velocity distributions of fluid particles in turbulent flows and the stochastic behavior of exchange rates [4],[5].

Let's $v(x,t)$ is a two variables function and assume that it can be demonstrated as a product of two functions which are single variable $v(x,t) = y(x)z(t)$. By making use of differential transform properties, $v(x,t)$ can be written as

$$v(x,t) = \sum_{i=0}^{\infty} Y(i) x^i \sum_{j=0}^{\infty} Z(j) t^j = \sum_{k=0}^{\infty} V_k(x) t^k \tag{3}$$

where $V_k(x)$ called $t$-dimensional spectrum function of $v(x,t)$ [1-3]. If the function $v(x,t)$ is analytic and differentiable continuously with respect to time $t$ and space $x$, so we can displayed $V_k(x)$ as

$$V_k(x) = \frac{1}{k!} \left[ \frac{\partial^k}{\partial t^k} v(x,t) \right]_{t=0} \tag{4}$$

Here, the lowercase $v(x,t)$ represents the original function and the uppercase $V_k(x)$ stand for the transformed.

From (4), we define the differential inverse transform of $V_k(x)$

$$v(x,t) = \sum_{k=0}^{\infty} V_k(x) t^k \qquad (5)$$

and then to compose (4) and (5), we get the numerical solution of $v(x,t)$ as below

$$v(x,t) = \sum_{k=0}^{\infty} \frac{1}{k!} \left[ \frac{\partial^k}{\partial t^k} v(x,t) \right]_{t=0} t^k. \qquad (6)$$

In the way of our utilization for the rest of the paper, from [1-3], the basic mathematical theorems taken by RDTM can be obtained as follows and their proofs are showed in [18].

**Theorem 1:** If $u(x,t)$ is two variables function, then transformation form is

$$U_k(x) = \frac{1}{k!} \left[ \frac{\partial^k}{\partial t^k} u(x,t) \right]_{t=0}.$$

**Theorem 2:** If $w(x,t) = u(x,t) \pm v(x,t)$, then transformation form $W_k(x) = U_k(x) \pm V_k(x)$.

**Theorem 3:** If $w(x,t) = \alpha u(x,t)$, then transformation form $W_k(x) = \alpha U_k(x)$, $\alpha$ constant.

**Theorem 4:** If $w(x,t) = x^m t^n$, then transformation form $W_k(x) = x^m \delta(k-n), \delta(k) = \begin{cases} 1, k=0 \\ 0, k \neq 0 \end{cases}$.

**Theorem 5:** If $w(x,t) = x^m t^n u(x,t)$, then transformation form $W_k(x) = x^m U_{k-n}(x)$.

**Theorem 6:** If $w(x,t) = u(x,t) v(x,t)$, then transformation form

$$W_k(x) = \sum_{r=0}^{k} U_r(x) V_{k-r}(x) = \sum_{r=0}^{k} V_r(x) U_{k-r}(x).$$

**Theorem 7:** If $w(x,t) = \frac{\partial^r}{\partial t^r} u(x,t)$, then transformation form $W_k(x) = (k+1)\ldots(k+r) U_{k+r}(x)$.

**Theorem 8:** If $w(x,t) = \frac{\partial}{\partial x} u(x,t)$, then transformation form $W_k(x) = \frac{d}{dx} U_k(x)$.

**Theorem 9:** If $w(x,t) = \frac{\partial^2}{\partial x^2} u(x,t)$, then transformation form $W_k(x) = \frac{d^2}{dx^2} U_k(x)$.

## 2. Numerical Considerations

In this part of paper, three illustrative examples which are wave-like equation with variable coefficients as in type of (1)-(2) are solved by RDTM for showing errors, convergence and efficient solutions.

**2.1. Example:** Let's at first, consider the nonlinear two dimensional wave-like equations which is a form of (1) with variable coefficients [4], [5],

$$v(x,y,t)_{tt} = \frac{\partial^2}{\partial x \partial y}\left(v(x,y,t)_{xx} v(x,y,t)_{yy}\right)$$

$$-\frac{\partial^2}{\partial x \partial y}\left(xyv(x,y,t)_x v(x,y,t)_y\right) - v(x,y,t) \quad (7)$$

with initial conditions

$$v(x,y,0) = e^{xy} \text{ and } v(x,y,0)_t = e^{xy}. \quad (8)$$

If we rearrange the equation (7), it can be written as following

$$v_{tt} = v_{xxxy}v_{yy} + v_{xxy}v_{yyx} + v_{xxx}v_{yyy} + v_{xx}v_{yyyx} - v_x v_y - xv_{xx}v_y - xv_x v_{xy}$$
$$- yv_{xy}v_y - xyv_{xxy}v_y - xy(v_{xy})^2 - yv_x v_{yy} - xyv_{xx}v_{yy} - xyv_x v_{yyx} - v \quad (9)$$

By applying the RDTM process for equation (9), we get reduced transformation form as below

$$(k+1)(k+2)V_{k+2}(x,y) = \sum_{r=0}^{k} \frac{d^2}{dy^2} V_r(x,y) \frac{d^4}{dx^3 dy} V_{k-r}(x,y)$$

$$+ \sum_{r=0}^{k} \frac{d^3}{dx^2 dy} V_r(x,y) \frac{d^3}{dy^2 dx} V_{k-r}(x,y) + \sum_{r=0}^{k} \frac{d^3}{dx^3} V_r(x,y) \frac{d^3}{dy^3} V_{k-r}(x,y)$$

$$+ \sum_{r=0}^{k} \frac{d^2}{dx^2} V_r(x,y) \frac{d^4}{dy^3 dx} V_{k-r}(x,y) - \sum_{r=0}^{k} \frac{d}{dx} V_r(x,y) \frac{d}{dy} V_{k-r}(x,y) \quad (10)$$

$$- x \sum_{r=0}^{k} \frac{d}{dy} V_r(x,y) \frac{d^2}{dx^2} V_{k-r}(x,y) - x \sum_{r=0}^{k} \frac{d}{dx} V_r(x,y) \frac{d^2}{dxdy} V_{k-r}(x,y)$$

$$- y \sum_{r=0}^{k} \frac{d}{dy} V_r(x,y) \frac{d^2}{dxdy} V_{k-r}(x,y) - xy \sum_{r=0}^{k} \frac{d}{dy} V_r(x,y) \frac{d^3}{dx^2 dy} V_{k-r}(x,y)$$

$$- xy \sum_{r=0}^{k} \frac{d^2}{dxdy} V_r(x,y) \frac{d^2}{dxdy} V_{k-r}(x,y) - y \sum_{r=0}^{k} \frac{d}{dx} V_r(x,y) \frac{d^2}{dy^2} V_{k-r}(x,y)$$

$$- xy \sum_{r=0}^{k} \frac{d^2}{dx^2} V_r(x,y) \frac{d^2}{dy^2} V_{k-r}(x,y) - xy \sum_{r=0}^{k} \frac{d}{dx} V_r(x,y) \frac{d^3}{dy^2 dx} V_{k-r}(x,y) - V_k(x,y)$$

and also the initial conditions of equation (8) are transformed

$$V_0(x,y) = e^{xy}, V_1(x,y) = e^{xy} \quad (11)$$

By substituting initial conditions into (10) and therewithal we apply the RDTM process as in the (4)-(6), then approximate solution of $v(x,y,t)$ is obtained following

$$v(x,y,t) = \sum_{k=0}^{\infty} V_k(x,y) t^k = e^{xy}\left(1 + t - \frac{t^2}{2} - \frac{t^3}{3!} + \frac{t^4}{4!} + \frac{t^5}{5!} - \frac{t^6}{6!} - \frac{t^7}{7!} + \ldots\right) \quad (12)$$
$$\cong e^{xy}(\sin(t) + \cos(t))$$

Thus, the result of (12) converges the exact solution of (7).

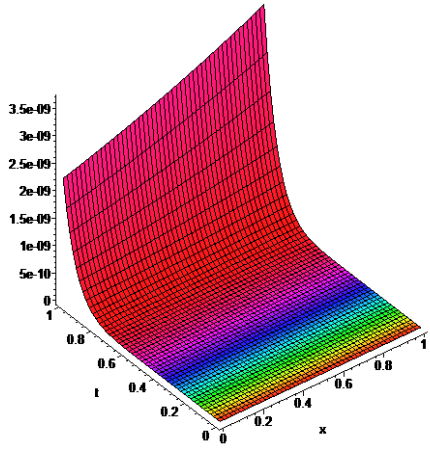

Figure 1: Absolute errors of example 2.1 between six terms RDTM solution and exact solution at y=0.5

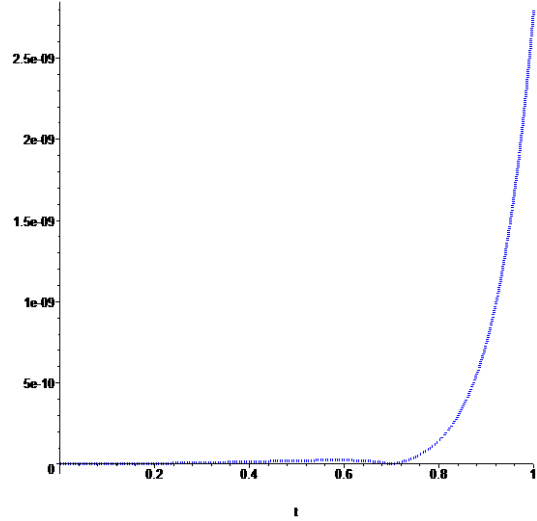

Figure 2: For y=x=0.5, absolute errors of example 2.1 between six terms RDTM solution and exact solution

Table 1: Absolute errors of four terms approximations of RDTM for example 2.1 which is better than to compare with Ref. [4] and Ref. [5].

| t/x,y | 0.1 | 0.2 | 0.3 | 0.4 | 0.5 | 0.6 | 0.7 | 0.8 | 0.9 | 1.0 |
|---|---|---|---|---|---|---|---|---|---|---|
| 0.1 | 0,25E-12 | 0,26E-12 | 0,27E-12 | 0,29E-12 | 0,31E-12 | 0,35E-12 | 0,4E-12 | 0,46E-12 | 0,56E-12 | 0,68E-12 |
| 0.2 | 0,6553E-10 | 0,6752E-10 | 0,7098E-10 | 0,7613E-10 | 0,833E-10 | 0,9298E-10 | 0,1059E-9 | 0,12304E-9 | 0,14582E-9 | 0,17633E-9 |
| 0.3 | 0,16967E-8 | 0,17484E-8 | 0,1838E-8 | 0,19713E-8 | 0,21569E-8 | 0,24077E-8 | 0,2742E-8 | 0,31857E-8 | 0,3776E-8 | 0,45662E-8 |
| 0.4 | 0,17117E-7 | 0,17638E-7 | 0,18542E-7 | 0,19887E-7 | 0,2176E-7 | 0,2429E-7 | 0,27662E-7 | 0,32139E-7 | 0,38094E-7 | 0,46065E-7 |
| 0.5 | 0,10301E-6 | 0,10614E-6 | 0,11159E-6 | 0,11968E-6 | 0,13095E-6 | 0,14617E-6 | 0,16647E-6 | 0,19341E-6 | 0,22925E-6 | 0,27722E-6 |
| 0.6 | 0,44704E-6 | 0,46065E-6 | 0,48427E-6 | 0,51938E-6 | 0,5683E-6 | 0,63438E-6 | 0,72245E-6 | 0,83936E-6 | 0,9949E-6 | 0,12031E-5 |
| 0.7 | 0,15481E-5 | 0,15953E-5 | 0,16771E-5 | 0,17987E-5 | 0,19681E-5 | 0,21969E-5 | 0,25019E-5 | 0,29068E-5 | 0,34454E-5 | 0,41664E-5 |
| 0.8 | 0,45445E-5 | 0,4682 E-5 | 0,4923 E-5 | 0,528 E-5 | 0,57772E-5 | 0,6449 E-5 | 0,73443E-5 | 0,85328E-5 | 0,10114E-4 | 0,1223 E-4 |
| 0.9 | 0,11758E-4 | 0,12116E-4 | 0,12737E-4 | 0,1366E-4 | 0,14947E-4 | 0,16685E-4 | 0,19001E-4 | 0,22076E-4 | 0,26166E-4 | 0,31642E-4 |
| 1.0 | 0,27533E-4 | 0,28371E-4 | 0,29826E-4 | 0,31989E-4 | 0,35001E-4 | 0,39071E-4 | 0,44495E-4 | 0,51696E-4 | 0,61276E-4 | 0,74097E-4 |

**2.2. Example:** Secondly, we take into account of the nonlinear wave-like equations which is a form of (1) with variable coefficients [4],[5] as follows

$$v(x,t)_{tt} = v(x,t)^2 \frac{\partial^2}{\partial x^2}\left(v(x,t)_x v(x,t)_{xx} v(x,t)_{xxx}\right)$$
$$+ \left(\frac{\partial v(x,t)}{\partial x}\right)^2 \frac{\partial^2}{\partial x^2}\left(\left(\frac{\partial^2 v(x,t)}{\partial x^2}\right)^3\right) - 18v(x,t)^5 + v(x,t) \quad (13)$$

with initial conditions

$$v(x,0) = e^x \text{ and } v(x,0)_t = e^x. \quad (14)$$

Again, if we apply all the (7)-(9) operations for the equation (13), it can be written as below

$$v_{tt} = v^2 \left(3v_{xx}(v_{xxx})^2 + 2(v_{xx})^2 v_{xxxx} + 3v_x v_{xxx} v_{xxxx} + v_x v_{xx} v_{xxxxx}\right)$$
$$+ (v_x)^2 \left(6v_{xx}(v_{xxx})^2 + 3(v_{xx})^2 v_{xxxx}\right) - 18v^5 + v \quad (15)$$

Hence, by using reduced differential transformation like (4)-(6) for on the equation (15), we write down transformed form

$$(k+1)(k+2)V_{k+2}(x) = 3\sum_{r=0}^{k}\sum_{s=0}^{k-r}\sum_{m=0}^{k-r-s}\sum_{n=0}^{k-r-s-m} V_r(x)V_s(x)\frac{d^2}{dx^2}V_m(x)\frac{d^3}{dx^3}V_n(x)\frac{d^3}{dx^3}V_{k-r-s-m-n}(x)$$
$$+ 2\sum_{r=0}^{k}\sum_{s=0}^{k-r}\sum_{m=0}^{k-r-s}\sum_{n=0}^{k-r-s-m} V_r(x)V_s(x)\frac{d^2}{dx^2}V_m(x)\frac{d^2}{dx^2}V_n(x)\frac{d^4}{dx^4}V_{k-r-s-m-n}(x)$$
$$+ 3\sum_{r=0}^{k}\sum_{s=0}^{k-r}\sum_{m=0}^{k-r-s}\sum_{n=0}^{k-r-s-m} V_r(x)V_s(x)\frac{d}{dx}V_m(x)\frac{d^3}{dx^3}V_n(x)\frac{d^4}{dx^4}V_{k-r-s-m-n}(x)$$
$$+ \sum_{r=0}^{k}\sum_{s=0}^{k-r}\sum_{m=0}^{k-r-s}\sum_{n=0}^{k-r-s-m} V_r(x)V_s(x)\frac{d}{dx}V_m(x)\frac{d^2}{dx^2}V_n(x)\frac{d^5}{dx^5}V_{k-r-s-m-n}(x) \quad (16)$$
$$+ 6\sum_{r=0}^{k}\sum_{s=0}^{k-r}\sum_{m=0}^{k-r-s}\sum_{n=0}^{k-r-s-m} \frac{d}{dx}V_r(x)\frac{d}{dx}V_s(x)\frac{d^2}{dx^2}V_m(x)\frac{d^3}{dx^3}V_n(x)\frac{d^3}{dx^3}V_{k-r-s-m-n}(x)$$
$$+ 3\sum_{r=0}^{k}\sum_{s=0}^{k-r}\sum_{m=0}^{k-r-s}\sum_{n=0}^{k-r-s-m} \frac{d}{dx}V_r(x)\frac{d}{dx}V_s(x)\frac{d^2}{dx^2}V_m(x)\frac{d^2}{dx^2}V_n(x)\frac{d^4}{dx^4}V_{k-r-s-m-n}(x)$$
$$- 18\sum_{r=0}^{k}\sum_{s=0}^{k-r}\sum_{m=0}^{k-r-s}\sum_{n=0}^{k-r-s-m} V_r(x)V_s(x)V_m(x)V_n(x)V_{k-r-s-m-n}(x) + V_k(x)$$

and initial conditions transform

$$V_0(x) = e^x, V_1(x) = e^x. \quad (17)$$

From reduced differential inverse transform process of (4)-(6), we obtain the numerical solution of (13)

$$v(x,t) = \sum_{k=0}^{\infty} V_k(x) t^k = v_0(x,t) + v_1(x,t) + v_2(x,t) + v_3(x,t) + \cdots$$
$$= e^x \left(1 + t + \frac{t^2}{2!} + \frac{t^3}{3!} + \frac{t^4}{4!} + \frac{t^5}{5!} + \frac{t^6}{6!} + \frac{t^7}{7!} + \frac{t^8}{8!} + \frac{t^9}{9!} \cdots \right) \cong e^{x+t} \qquad (18)$$

which converges efficiently to exact solution.

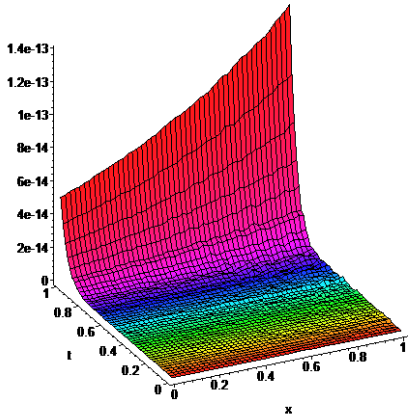
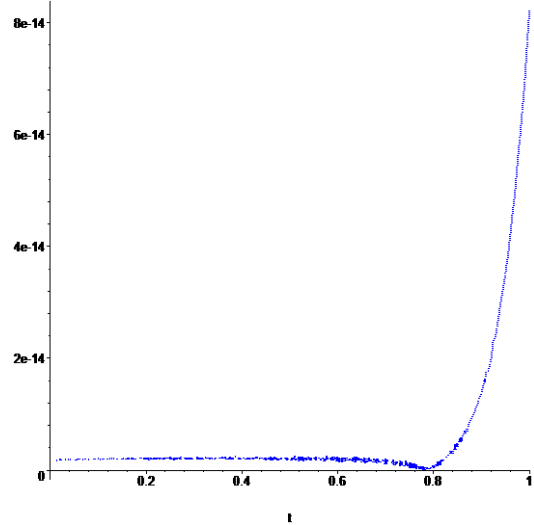

**Figure 3:** Absolute errors of example 2.2 between eight terms RDTM solution and exact solution

**Figure 4:** For x=0.5, absolute errors of example 2.2 between eight terms RDTM solution and exact solution

**Table 2:** Absolute errors of eight terms approximations of RDTM for example 2.2 which is better than to compare with Ref. [4] and Ref. [5].

| t/x | 0.2 | 0.4 | 0.6 | 0.8 | 1.0 |
|-----|-----|-----|-----|-----|-----|
| **0.2** | 0,2E-13 | 0,2E-13 | 0,2E-13 | 0,2E-13 | 0,3E-13 |
| **0.4** | 0 | 0 | 0 | 0,1E-13 | 0,1E-13 |
| **0.6** | 0,1E-13 | 0,1E-13 | 0,1E-13 | 0 | 0 |
| **0.8** | 0,2E-13 | 0,3E-13 | 0,2E-13 | 0,3E-13 | 0,4E-13 |
| **1.0** | 0,6E-13 | 0,7E-13 | 0,9E-13 | 0,11E-12 | 0,12E-12 |

**2.3. Example:** And finally, we handle the nonlinear wave-like equations which is a form of (1) with variable coefficients [4], [5]

$$v(x,t)_{tt} = x^2 \frac{\partial}{\partial x}\left(v(x,t)_x \, v(x,t)_{xx}\right) - x^2 \left(\frac{\partial^2 v(x,t)}{\partial x^2}\right)^2 - v(x,t) \qquad (19)$$

with initial conditions

$$v(x,0) = 0 \text{ and } v(x,0)_t = x^2. \tag{20}$$

We rewrite the equation (19) under cover of the example 2.1 and 2.2 to as

$$v_{tt} = x^2\left((v_{xx})^2 + v_x v_{xxx}\right) - x^2(v_{xx})^2 - v \tag{21}$$

We recall that the processing steps of (7)-(9), (13)-(15) and also apply the reduced differential transform to equation (21), transformed form is obtained easily

$$(k+1)(k+2)V_{k+2}(x) = x^2\left(\sum_{r=0}^{k}\frac{d^2}{dx^2}V_r(x)\frac{d^2}{dx^2}V_{k-r}(x) + \sum_{r=0}^{k}\frac{d}{dx}V_r(x)\frac{d^3}{dx^3}V_{k-r}(x)\right)$$
$$- x^2\sum_{r=0}^{k}\frac{d^2}{dx^2}V_r(x)\frac{d^2}{dx^2}V_{k-r}(x) - V_k(x) \tag{22}$$

and initial conditions transform

$$V_0(x) = 0, V_1(x) = x^2. \tag{23}$$

The differential inverse transform of $V_k(x)$ gives an approximate solution as RDTM operations in (4)-(6)

$$v(x,t) = \sum_{k=0}^{\infty}V_k(x)t^k = v_0(x,t) + v_1(x,t) + v_2(x,t) + v_3(x,t) + \cdots$$
$$= x^2\left(t - \frac{t^3}{3!} + \frac{t^5}{5!} - \frac{t^7}{7!} + \ldots\right) \cong x^2 \sin t \tag{24}$$

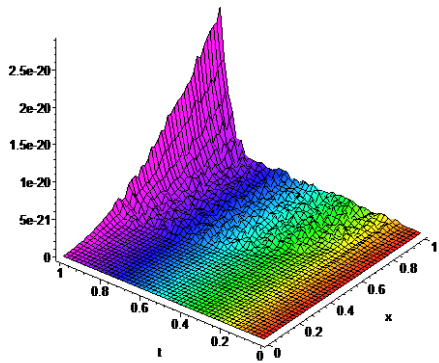

**Figure 5:** Absolute errors of example 2.3 between ten terms RDTM solution and exact solution

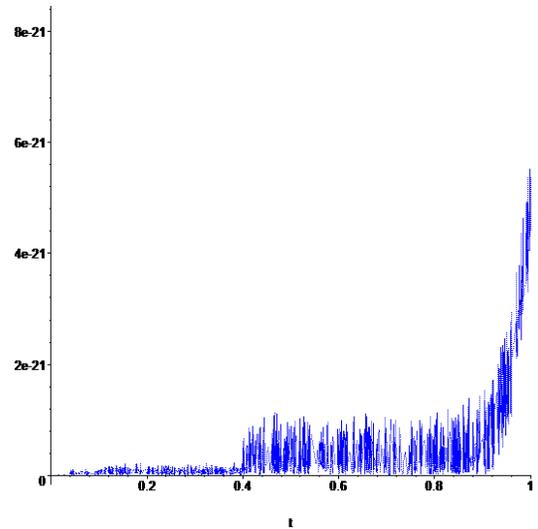

**Figure 6:** For x=0.5, absolute errors of example 2.3 between ten terms RDTM solution and exact solution

Table 3: Absolute errors of ten terms approximation of RDTM for example 2.3 which is better than to compare with Ref. [4] and Ref. [5].

| t/x | 0.2 | 0.4 | 0.6 | 0.8 | 1.0 |
|-----|-----|-----|-----|-----|-----|
| 0.2 | 0,2E-26 | 0 | 0,1E-25 | 0,1E-24 | 0,1E-24 |
| 0.4 | 0,1E-25 | 0,2E-25 | 0 | 0,1E-24 | 0,2E-24 |
| 0.6 | 0,3E-25 | 0,5E-25 | 0,1E-24 | 0,4E-24 | 0,3E-24 |
| 0.8 | 0,72E-23 | 0,286E-22 | 0,649E-22 | 0,1153E-21 | 0,1803E-21 |
| 1.0 | 0,78138E-21 | 0,31255E-20 | 0,70323E-20 | 0,125021E-19 | 0,195343E-19 |

## 3. Discussion

In this article, we applied the reduced differential transform method (RDTM), which has an advantage to provide an analytical approximation to the solution, usually an exact solution, in a rapidly convergent sequence, for nonlinear wave-like equations with variable coefficients. RDTM can be performed very easily and it is more effective and reliable, when it is compared with most famous techniques (Adomian Decomposition and Homotopy Perturbation) as in [4], [5]. Additionally, RDTM is faster than ADM-HPM to solve this type of equations.

So, our results show that the presented method is powerful technique and provides high accuracy for solving wave-like equations.